\documentstyle[12pt,amssymb,twoside]{article}
\newtheorem{pr}{Proposition}[section]
\newtheorem{th}[pr]{Theorem}
\newtheorem{co}[pr]{Corollary}
\newtheorem{lm}[pr]{Lemma}
\newtheorem{ex}[pr]{Example}
\newtheorem{df}[pr]{Definition}
\newtheorem{re}[pr]{Remark}

\begin{document}
\pagestyle{myheadings} \markboth{\footnotesize{M. Abolghasemi, A.
Rejali and H. R. E. Vishki  }}{\footnotesize{ On the
transformation semitopological semigroup...}}
\title{On the transformation semitopological
semigroup}
\author{\footnotesize{ M. Abolghasemi, A. Rejali and H. R. E.
Vishki }}
\date{}
\date{}
\maketitle 

\begin{abstract}In this paper we introduce the notion of weighted (weakly) almost periodic compactification
of a semitopological semigroup and generalize this notion to
corresponding notion for transformation semigroup.The inclusion
relation and equality of some well-known function spaces on a
weighted transformation semigroup is also investigated.
\\\\
\footnotesize\emph{\emph{2000} MSC}: 22A15, 32J05.
\end{abstract}
\addtolength{\baselineskip}{.4mm}
\maketitle

\section*{Introduction}
\quad ~Let $\omega$ be a weight on a semitopological semigroup $X$
with identity $e$ and $S$ be a semigroup. The pair $(S,X)$ is
semitopological transformation semigroup and $AP(X)\\ (WAP(X))$ is
denoted almost periodic (weakly almost periodic) function space.
For a topological group $G.$

Burckel proved that $G$ is compact if and only if bounded
continuous functions on $G$ are weakly almost periodic
\cite{Burckel}. Granirer improved the result to $G$ is compact if
and only if all the bounded uniformly continuous functions on $G$
are weakly almost periodic \cite{Granirer}. Dzinotyiweyi proved a
same result for very large class of topological semigroup. H. D.
Junghenn generalize the notion of (weakly) almost periodic
compactificatin of semitopological semigroup to the corresponding
notion for transformation semigroup and show that if $S$ has a
(weakly) almost periodic compactification then $(S,X)$ has a
(weakly) almost periodic compactification \cite{Junghenn}.

A. Rejali gave a modification of definition of weighted (weakly)
almost periodic function spaces and proved that
$AP(G,\omega)=C(G,\omega)$ when $G$ is compact group
\cite{Rejali}.

The organization of this paper is as follows.\\
In section 1 we introduce preliminaries and some definition.\\
In section 2 we generalize definition of (weakly) almost periodic
function space in \cite{Rejali} and introduce weighted
transformation functions spaces on weighted transformation
semigroup and prove that each function in $C(X,\omega)$ is almost
periodic if $S$ and $X$ are compact.\\Also we show that whenever a
member of $C(X,\omega)$ or $WAP(X,\omega)$ belong to
$AP(X,\omega)$. We will give examples of the weighted
semitopological transformation semigroup and study the
equicontiniuty  of this function spaces.

In section 3 we introduce a notion of weighted transformation
compactification and finally we give the main result of the paper
which states; each weakly almost periodic function on $X$ with
respect to $S$ is left norm continuous.

\section{Preliminaries  }
 \quad~ Let $(S,X,\omega)$ be a
weighted transformation semigroup in which $S$ is a semigroup and
$\omega$ is a weight on $X$.($\omega :X \longrightarrow
\mathbb{R}^+ $ is called a weight on $X$ if $\omega (xy)\leq
\omega(x) \omega(y) ~~x,y\in X$ and $\omega , \omega ^{-1}$ are
bounded on compact subsets of $X$).

A map $(s,x)\longrightarrow sx : S \times X \longrightarrow X$ is
called the action of S on X if $s(tx)=(st)x.$ For all $s,t \in
S~~, x\in X$.

If $T$ is a subsemigroup of $S$ and $Y$ is a $T$- invariant subset
of $X$, i.e. \\
$TY=\{ty :t\in T , y\in Y\}\subseteq Y$. Then $(T,Y,\omega)$ is
called a sub-transformation semigroup of $(S,X,\omega).$

For topological spaces $S$ and $X$ , $(T,Y,\omega)$ is dense in
$(S,X,\omega)$  if $T$ and $Y$ are dense in $S$ and $X$
respectively. $C(X)$ is the $C^*$-algebra of bounded continuous
real-valued function on $X$ with supnorm.

If $\mathcal{F}$ is a $C^*$-subalgebra of $C(X)$ then
$X^{\mathcal{F}}$ denotes the spectrum of $\mathcal{F}$
 (the set of all multiplicative means on$\mathcal{F}$) which is $\omega^*$ -compact in $\mathcal{F}^*$\\

The evaluation map $\epsilon :X\longrightarrow X^{\mathcal{F}}$
with $\omega^*$-dense image is defined by $\epsilon(x)(f)=f(x)$
for $ x\in X , \ f\in \mathcal{F}.$\\

For $s,t \in S$ and  $ x\in X$, define $\lambda _
s:S\longrightarrow S $ and $\rho _ t:S\longrightarrow S $ by
$\lambda _ s(t)=st=\rho _ t(s)$.Also $\hat\lambda _s
:X\longrightarrow X$ and $\hat\rho _x :S\longrightarrow X$ by
$\hat\lambda _s (x)= \hat\rho _x (s)=sx.$\\ These maps are called
 translation maps.

 $(S,X,\omega)$ is left topological if $S$ and $X$ are topological space and $\lambda _s , \hat\lambda _s $
  are continuous for all $s\in S$, right topological if  $\rho _s $ and $\hat\rho _x $ are continuous
  for all $s \in S , x \in X$, and semitopological if it is both left and right
  topological.

 A semitopological semigroup $S$ is topologically
 right (left) simple if for each $s \in S; ~sS ~(Ss)$ is dense in  $S$.\\
   $(S,X,\omega) $ is called topological if the multiplication in S and the action of S on X are
   continuous.

    If there is a separately (jointly) continuous action of a semitopological (topological)
    group $S$ with identity $e$ on a topological space $X$ and $ex=x$  for all $x\in
    X$,
    then ,$(S,X,\omega)$ is called a semitopological(resp.topological) transformation group.

$(S,X,\omega)$ is said to be compact if so are S and X. We put ,\\
 $L_s =(\lambda _s)^*$ and $R_s =(\rho _s)^*$ on C(S).
  $\hat L _s=(\hat \lambda _s )^* $ and $\hat R_x =(\hat \rho
  _x)^*$ on C(X)\\ where $(\lambda _s)^* $is dual of $(\lambda _s)$ and $(\rho
  _s)^*$ is dual of $\rho _s$.   clearly $\hat L_{st} =\hat L_s
   \hat L_t$ and $\hat R_{sx}=R_s \hat R _x$  and  $L _s \hat R _x=\hat R_x \hat L
   _s$

   A subset $\mathcal{F}$ of C(S) is called translation invariant
   if $L_s \mathcal{F} \bigcup$ $R_s \mathcal{F} \subseteq
   \mathcal{F}$  for $s \in S.$

   Let $A=\{f_\alpha :\alpha \in I\} $
   be a family of real functions on topological space $X$. $A$ is
   called equcontinuous if for every $\epsilon >0$ and $x_0 \in X$ there
   exit a neighborhood $N$ of $x_0$ such that $|f_\alpha (x)-f_\alpha
   (x_0)|<\varepsilon$,
   for every $x\in N$ and $\alpha \in I.$\\

\newpage
\section{ Weighted Function Spaces }

\begin{df}\emph{ Let $(S,X,\omega)$ be a weighted semitopological transformation
semigroup, and $e$ be the identity of X. Then:\\
 for weight $\omega $, put $~~~~~\Omega :X \times X  \longrightarrow (0,1]$ defined $\Omega(x,y)=\frac{\omega(xy)}{\omega(x)\omega(y)} $
and $_s\Omega:X \longrightarrow (0,1]$ by $_s\Omega(x)=
\Omega(se_x,x )$ and $\Omega_x:S \longrightarrow (0,1]$ by
$\Omega_x(s)=\Omega(se_x,x ).$}
\end{df}

\textbf{Notation:}  Always assume that $\Omega$ is separately
continuous and we will use from $_sf, f_x$ instead of $\hat L_sf$
and $\hat R_x$.

Now we introduce some weight function spaces.
\[C(X,\omega)=\{f\in C(X) : \frac{f}{\omega} \in C(X)\}\] and for
every $f\in C(X, \omega)$ define:
\[\|f\|_\omega = sup\{\frac{|f(x)|}{\omega(x)}:  x\in
X\}.\]\\
$L_{S_\omega}f=\{_s(\frac{f}{\omega}) _s\Omega  : s\in S \}$ and
$R_{S_\omega}f=\{(\frac{f}{\omega})_x \Omega_x  : x\in X \}$ \\
 $AP(X,\omega)=\{f \in C(X,\omega) :
\{\hat{L}_s(\frac{f}{\omega})_s\Omega: s \in S \}$ is relatively
compact in norm topology on $C(X)$ $\}$ \\
$WAP(X,\omega)=\{f \in C(X,\omega) :
\{\hat{L}_s(\frac{f}{\omega})_s\Omega: s \in S \}$ is relatively
compact in weak topology on $C(X)$\}$ \\
LUC (X,\omega)=\{f \in C(X,\omega)$ :The map $s\longrightarrow
\hat{L}_s(\frac{f}{\omega}) _s\Omega :
S\longrightarrow C(X)$ is norm continuous\}\\
$RUC(X,\omega)=\{f \in C(X,\omega)$ :The map $x\longrightarrow
\hat{R}_x(\frac{f}{\omega}) \Omega_x : X\longrightarrow C(S)$ is
norm continuous$\}$\\
\begin{pr}Let $f\in C(X,\omega)$. Then
$f \in AP(X,\omega)$  ($WAP(X,\omega)$ ) if\\
$\{\frac{_sf}{\omega (se)} : s \in S \}$is relatively (weakly)
compact in $C(X,\omega).$
\end{pr}
\textbf{Proof.}\\ Let  $x\in X$. Then :\\
\[_s(\frac{f}{\omega})_s\Omega(x)=\frac{f(sx)}{\omega(sx)}
\Omega(se,x)=\frac{f(sx)}{\omega(sx)}
\frac{\omega((se)x)}{\omega(se)\omega(x)}=\frac{f(sx)}{
\omega(se)\omega(x)}=(\frac{_sf}{\omega(se)\omega})(x)\]\\
therefore:\\
\[_s(\frac{f}{\omega})_s\Omega=\frac{_sf}{\omega(se) \omega}\]and so if $f \in C(X, \omega) $ then: \\
$ f \in AP(X,\omega)$  (resp$_.  WAP(X,\omega))$ if
$\{\frac{_sf}{\omega(se)} : s \in S\}$ is relatively (weak)compact
in $C(X,\omega).$\\

\begin{ex} \emph{If $S$ be a right zero semigroup and $S=X$ then $AP(X,\omega)=C(X,\omega)$.}
\end{ex}
\begin{ex} \emph{let $(S,X,\omega)$ be a weighted semitopological transformation
semigroup and $\omega$ be a multiplicative weight on $X$ i.e
$(\omega (xy)= \omega(x) \omega(y),x,y \in X)$. Then
$AP(X,\omega)=AP(X)$ and $WAP(X,\omega)=WAP(X).$}
\end{ex}

\begin{re}
\emph{All of these function spaces are translation invariant
$C^*$- subalgebras of $C(X, \omega)$ containing the constant
function and clearly $AP(X, \omega) \subseteq WAP(X, \omega).$}
\end{re}
\begin{lm}
Let $(S,X,\omega)$ be a weighted semitopological transformation
semigroup. Then:\\
 \textbf{1)} If $(S,X,\omega)$ is compact, then $WAP(X,\omega)=C(X,\omega).$\\
 \textbf{2)} If $(S,X,\omega)$ is a compact, then
$AP(X,\omega)=C(X,\omega).$\\
 \textbf{3)} If $f \in WAP(X,\omega)$, then the map
$ s\longrightarrow _s(\frac{f}{\omega}) _s\Omega :S\longrightarrow
WAP(X,\omega)$ is continuous in the weak topology.\\
\textbf{4)} If $f \in WAP(X,\omega)$, then the map $
s\longrightarrow _s(\frac{f}{\omega})
  _s\Omega :S\longrightarrow AP(X,\omega)$ is continuous in the norm
  topology.

\emph{we prove that only 1 and 3. The proofs of 2 and 4 are similar.}\\
\textbf{Proof.}\\
 \textbf{1)} \emph{Let $(S,X,\omega)$ be compact and $f\in
C(X,\omega)$,since the map \\
$ s\longrightarrow _s(\frac{f}{\omega}) _s\Omega $ is pointwise
continuous , $\{ _s(\frac{f}{\omega}) _s\Omega :s\in S \}$ is
compact in the pointwise topology of C(X) ,and since this topology
agrees with the weak topology on norm bounded pointwise compact of
$C(X,\omega)$ (\cite{Junghenn}). $
 \{_s(\frac{f}{\omega}) _s\Omega :s\in S\}$ is relatively  weakly
compact. therefore $f\in WAP(X,\omega).$\\}
 \textbf{3)}\emph{ Let
$(s_\alpha)$ a net in S and converges to s. since $f \in
 WAP(X,\omega)$ then $\{_s(\frac{f}{\omega}) _s\Omega :s\in
S\}$ is relatively weakly compact therefore $ (s_\alpha
(\frac{f}{\omega})s_\alpha \Omega )$ has unique weak limit point
in $C(X)$. Therefore $ _{s_\alpha} (\frac{f}{\omega}) _{s_\alpha}
\Omega \longrightarrow _s(\frac{f}{\omega}) _s \Omega$\\}

\end{lm}
 \begin{pr} Let $(S,X.\omega)$ be a weighted semitopological transformation
semigroup. Then:\\
\textbf{1)} $AP(X,\omega)\subseteq LUC(X,\omega) \cap RUC(X,\omega)$\\
\textbf{ 2)} If $S$ ($X$) is compact, then
$AP(X,\omega)=LUC(X,\omega)$ ( $RUC(X,\omega)$.\\
\textbf{3)} If$(S,X,\omega)$ is compact, then
 $AP(X,\omega)=LUC(X,\omega)=RUC(X,\omega)$.\\
\textbf{4)}If $(S,X,\omega)$ is compact and topological then
   $AP(X,\omega)=LUC(X,\omega)=RUC(X,\omega)=WAP(X,\omega)=C(X,\omega).$\\
\textbf{ 5)} If $S$ ($X$) is compact and Hausdorff, then
 \[C(X,\omega)= RUC(X,\omega)\ (~LUC(X,\omega))\] if and only if the
action of S on X is (jointly) continuous .\\
 \textbf{6)} If $(S,X,\omega)$ is compact and Hausdorff, then $LUC(
 X,\omega)=RUC(X,\omega)=C(X,\omega)$ if and only if the action of  S on X
 is continuous.\\
\textbf{Proof.}\\
\textbf{1)} \emph{ If $f\in AP(X,\omega)$ then the map
$s\longrightarrow _s(\frac{f}{\omega}) _s\Omega $ is continuous by
\cite{Dzinotyiweyi}. Therefore $f \in LUC(X,\omega)$.\\
Now if $ f\in RUC(X,\omega)$ the map $ x\longrightarrow
(\frac{f}{\omega})_x \Omega_x $ is continuous. Let $ U(x)=
(\frac{f}{\omega})_x \Omega_x $ and $ V(s)= _s
(\frac{f}{\omega})_s\Omega ~~~ $ for $ x\in X ,s\in S. $ Obviously
$(V(s))(x) =(U(x))(s).$ Therefore $ f\in RUC(X,\omega),$ hence$
\\AP(X,\omega)=LUC(X,\omega)\bigcap RUC(X,\omega).$\\}
\textbf{2)} \emph{Let $S $ be compact semigroup and $ f\in
LUC(X,\omega)$. Then the map $ s\longrightarrow
 _s(\frac{f}{\omega})_s\Omega $ is continuous.
 Therefore $\{_s(\frac{f}{\omega})_s\Omega : s\in S \} $ is compact
 (since $S$ is compact) so $ f\in AP(X,\omega)$.\\
 Note that if $X$ is compact, then the proof is similar.\\}
\textbf{3)}\emph{It follows from (2) and (1 lemme 2.6)\\}
 \textbf{4)}\emph{It follows from (3) and (2 lemme 2.6)\\}
 \textbf{5)}\emph{ we note that the action of S on X is continuous if and
 only if for $ f\in C(X,\omega)$ the map $(s,x)\longrightarrow f(sx)
 :S \times X\longrightarrow \mathbb{R} $
 is continuous (App. B 3 of \cite{Berg}). This equivalent to norm
 continuity the map $ x\longrightarrow
 (\frac{f}{\omega})_x \Omega_x$, when S is compact.\\}
\textbf{6)} \emph{If $ (S,X,\omega)$ is compact,Then by (1) lemma
2.6 and 3 proof is obviously.\\}

\end{pr}
\begin{pr}
Let $(S,X,\omega)$ be a weighted semitopological transformation
semigroup and $f\in C(X,\omega).$ Then:\\ \textbf{1)}$ f\in
RUC(X,\omega)$ if and only if $L_{S_\omega}f$ is
equicontinuous.\\
\textbf{2)}$f\in LUC(X,\omega)$ if and only if
$R_{S_\omega}f$ is equicontinuous.\\
\textbf{Proof.}\\
\emph{1) Let $f\in RUC(X,\omega)$, Then the map
  $x\longrightarrow (\frac{f}{\omega})_x \Omega_x : X\longrightarrow C(S)$ is
norm continuous. Hence for every $ \epsilon>0$ and $x_0 \in X $
there is a neghborhood $N$ of $x_0$ such that for every $x\in X;
~~~~ \|(\frac{f}{\omega})_x \Omega_x-(\frac{f}{\omega})_{x_0}
\Omega_{x_0}\|<\epsilon$. Hence for every $s\in S~~~ \mid
\frac{f}{\omega}(sx)\Omega_x(s)-\frac{f}{\omega}(sx_0)\Omega_{x_0}(s)\mid<\epsilon$
and so $ \mid
_s{(\frac{f}{\omega})}_s\Omega(x)-_s{(\frac{f}{\omega})}_s\Omega(x_0)\mid<\epsilon$.
Hence $L_{S_\omega}f$ is equicontinuous. The Converse is clear and
the proof of 2 is simillar.}
\end{pr}
\begin{co}
\textbf{1)}If $S$ is compact then for every $f\in AP(X,\omega)$
the set $ R_{S_{\omega}}f$ is equicontinuous.\\\textbf{2)}If $X$
is compact then for every $f\in AP(X,\omega), L_{S_{\omega}}f$ is
aquicontinuous.\\
 \textbf{Proof.} By Proposition (2.7)and Proposition (2.8) is immediate.
\end{co}
\begin{df} \emph{A weighted homomorphism of a weighted semitopological transformation semigroup $ (S,X,\omega_1)$ into
a weighted semitopological transformation semigroup
$(T,Y,\omega_2)$ is a pair $ (\sigma ,\eta)$, where $\sigma
:S\longrightarrow T $ is a continuous homomorphism and $\eta
:X\longrightarrow Y$ is a continuous homomorphism such that $ \eta
(sx)=\sigma(s)\eta(x)$ for $x\in X$ and  $s\in S.$ and
$\omega_2\circ \eta =\omega_1$.}
\end{df}
\begin{re}

 \emph{Let $ (\sigma ,\eta)$ be a homomorphism of $(S,X,\omega_1)$
into $(T,Y,\omega_2)$. Let\\
  $ f\in C(Y,\omega_2) , s \in S$ and $x \in X. $
Then :\\\textbf{1)}~~~~~$_s(\frac{\eta^*
(f)}{\omega_1})(x)=(\eta^* (\frac{f}{\omega_2})(sx))$
=$(\frac{f}{\omega_2})o \eta
(sx)=(\frac{f}{\omega_2})(\sigma(s)\eta(x))=\\_{\sigma(s)}(\frac{f}{\omega_2})o\eta(x)
=\eta^* (_{\sigma(s)}(\frac{f}{\omega_1}))(x)$, and so
:\[_s(\frac{\eta^*(f)}{\omega_1})=\eta^*(_{\sigma(s)}(\frac{f}{\omega_2}))\]
 \textbf{2)} $${_{s}\Omega_{1}=\eta^{\ast}
 (_{\sigma_s}\Omega_2)}$$}
\end{re}
\begin{lm} Let $(\sigma,\eta)$ be a homomorphism of $(S,X,\omega_1)$
into $(T,Y,\omega_2)$. where $\omega_2 \circ \eta =\omega_1$.
Then:\[ \eta^* (WAP(Y,\omega_2))\subseteq WAP(X,\omega_1) \cap
\eta^* (C(Y,\omega_2))\] and equality is hold if $\eta (X)$ is
dense in Y and $\sigma (S)$ is dense in T.\\
 \textbf{Proof.}\\\emph{ Let $ f\in WAP(Y,\omega_2)$ and
 $A=\{{_{\sigma s}(\frac{f}{\omega_2})}{_{\sigma s}\Omega_2 :~s\in S} \}.$\\
 Then~
 $\eta^*(A)$ is relatively weakly compact, and by remark (2.12) $\eta^*(f)\in WAP(X,\omega_1)$. Since
 $\eta^*(f)\in \eta^*(C(Y,\omega_2))$. Therefore \[\eta^*(f)\in
 WAP(X,\omega_1) \cap \eta^*(C(Y,\omega_2)).\]\\}

\emph{Now if $\eta(X)$ and $\sigma(S)$ are dense in Y and T,
respectively, then $\eta$ is an isometry and by remark (2.12) ~~~
$(\eta^*)^{-1} (cl(B)) = cl(A)$ ~~~~~~~~  (1)
\[ \emph{\emph{where}} ~~B=\{_{\sigma(s)}(\frac{f}{\omega_1})_{\sigma(s)}\Omega_1 :s\in
S. \}\] Let $g\in WAP(X,\omega_1)\cap \eta^* (C(Y,\omega_2)$. Then
$(\eta^*)^{-1}(g)\in C(Y,\omega_2)$ and by remark (1.7) and (1)
$(\eta^*)^{-1} (g) \in WAP (Y,\omega_2)$ and $g=\eta^*
(\eta^*)^{-1} (g)\in WAP (Y,\omega_2)$. Hence equality is hold.}
\end{lm}
\section{Weighted Transformation Compactification}
\begin{df}
 \emph{ Let $(S,X,\omega)$ be a weighted transformation semitopological
semigroup and $(T,Y,\bar{\omega})$ be a compact Hausdorff right
topological weighted semigroup. Further $(\phi ,\psi):
 (S,X,\omega) \longrightarrow (T,Y,\bar{\omega})$ be a continuous
homomorphism such that $\phi(S)$ and $\psi(S)$ are dense
semitopological sub transformation semigroup of T and Y
respectively. Then  $((\phi ,\psi),(T,Y,\bar{\omega}))$ is called
weighted  transformation compactification of $(S,X,\omega).$\\}
 \emph{ If $((\phi ,\psi),(T,Y,\bar{\omega}))$ and $((\phi^{'}
,\psi^{'}),(T^{'},Y^{'},\omega^{'}))$ are two weighted
 transformation compactification of $(S,X,\omega)$
and $(\pi , \gamma):(T,Y,\bar{\omega})\longrightarrow
(T^{'},Y^{'},\omega^{'})$ is a continuous homomorphism of
$(T,Y,\bar{\omega})$ onto $(T^{'},Y^{'},\omega^{'}).$ we say that
$((\phi ,\psi),(T,Y,\bar{\omega}))$ is an extension of $((\phi^{'}
,\psi^{'}),(T^{'},Y^{'},\omega^{'}))$}

\emph{ Let $(S,X,\omega)$ be a weighted transformation
semitopological semigroup and P be a property of compactification
of  $((\phi ,\psi),(T,Y,\bar{\omega}))$.Then a P-compactification
of $(S,X,\omega)$ is a compactification of $(S,X,\omega)$ that has
the given property P.}

\emph{If a P-compactification of $(S,X,\omega)$ that is extension
of each other P-compactifica-\\tion of $(S,X,\omega)$ is called
universal P-compactification.}
 \end{df}
 \begin{re}
 \emph{If $((\phi ,\psi),(T,Y,\bar{\omega}))$ is a weighted transformation compactification of
 $(S,X,\omega)$, then $\psi^\ast
 (\bar{\omega})=\omega$.}
\end{re}
\begin{ex} \emph{For a  weighted transformation semitopological semigroup
 $(S,X,\omega)$ $((\epsilon ,\delta ) (S^{WAP(S)},X^{WAP(X,\omega)},\bar{\omega}) $ is a
 universal semitopological compactification of $(S,X,\omega)$,
 where $\epsilon$ and $\delta$ are evaluation map on $ C(S)$ and $
 C(X,\omega) $ respectively.}
\end{ex}
\begin{ex}\emph{ Let $(S,X,\omega)$ be a semitopological semigroup where $S=X$ and $e$ be the identity of
$X$. Further let $M$ be clousure of ${e}$ in $X$, then  $M$ is a
closed congruence subsemigroup of $X$ (5.4 of \cite{Hewitt} and
(1.24) of \cite{Berg}). Put $H=\frac{X}{M}$. Hence $H$ is a
Hausdorff topological semigroup. Suppose that
$\pi:X\longrightarrow H$ denotes canonical homomorphism. For each
$f\in C(X,\omega) $ define $g\in C(H,\bar{\omega)}$ by
$g(\pi(x))=f(x).$ Since $f$ is constant on $H$ Then $g$ is
well-defined and hence $g\in C(H,\bar{\omega}).$ Thus $\pi^*
:C(H,\bar{\omega})\longrightarrow C(X,\omega)$ is an isometric
isomorphism of $C(H,\bar{\omega})$ on $C(X,\omega).$\\
Therefore $\pi^*(R_{\pi(X)}g)=R_Sf$ and it follows
 $\pi^*(WAP(H,\bar{\omega}))=WAP(X,\omega)$ and if
$WAP(X,\omega)=C(X,\omega)$, then $
WAP(H,\bar{\omega})=C(H,\bar{\omega})$}

\end{ex}
\begin{lm} Let $(S,X,\omega)$ be a compact Hausdorff weighted
semitopological transformation semigroup and $ T= \{t \in S : tS=S
, tX=X\}$ be a dense subsemigroup of S. Then for each $f\in
C(X,\omega)$ the map $s \longrightarrow _s(\frac{f}{\omega
})_s\Omega :S \longrightarrow C(X)$ is norm continuous at each
point of $T$.\\
\textbf{Proof.}

\emph{Let f $\in C(X,\omega)$. By (App. B1 of \cite{Berg}); there
is $s_0 \in S$ such that the function $(s,x) \longrightarrow
\frac{f}{\omega}(sx) :S \times X \longrightarrow \mathbb{R}$ is
jointly continuous at each point of $\{s_0\} \times X$. Further
By (App.B3 of \cite{Berg}) the set
$$N=\{s \in S :~
       \|_s (\frac{f}{\omega})- _{s_0}(\frac{f}{\omega})\|
< \frac{\epsilon}{2} \}$$  for every $\epsilon >0$ is a
neighborhood of $s_0$. }

\emph{Let $t_0\in T$ be arbitrary,then by definition of $T$ there
exist $u_0 \in S$ such that $t_0 u_0 =s_0$. Since $T$ is dense in
$S$ and $s_0 S=S$, then $t_0T$ is dense in $S$.\\Now chooses $t
\in T$ such that $t_0t \in N$ and put $V=\rho^{-1} _t (N)$. Then
$V$ is a neighborhood of $t_0$. For
every $s \in V$ we have :$\\
\|_s (\frac{f}{\omega})- _{t_0} (\frac{f}{\omega})\|=sup_{x\in X}
|\frac{f}{\omega}(sx)-\frac{f}{\omega}(t_0x)|=
sup_{x\in X}|\frac{f}{\omega}(stx)-\frac{f}{\omega}(t_0tx)|\leq \\
\|_{st} (\frac{f}{\omega})- _{t_0 u_0} (\frac{f}{\omega})\|+
\|_{t_0 u_0} (\frac{f}{\omega})- _{t_0 t} (\frac{f}{\omega})\|\leq
\frac{\epsilon}{2}+ \frac{\epsilon}{2} =\epsilon .$ since
$s_0=t_0u_0$ and $st , t_0t \in N .$\\
 Therefore the map $s
\longrightarrow _s(\frac{f}{\omega}) _s\Omega  :S \longrightarrow
C(X) $ is norm continuous at $t_0$.}
\end{lm}
\begin{co}\emph{If $(S,X,\omega)$ is a compact Hausdorff weighted semitopological transformation semigroup
and $T=\{ t \in S :tS=S ,tX=X \} $ be a dense subsemigroup of $S$,
then $(T,X,\omega)$ is a topological transformation subsemigroup
of $(S,X,\omega)$.}
\end{co}
\begin{df}\emph{ Let $(S,X,\omega)$ be weighted transformation
semigroup. Then $S$ acts surjectively (resp$_.$ topologically
surjectively ) on $X$ if for every $s \in S$ have $sX=X.$ (resp.
$_. sX$ is dense in $X$.)\\}
\end{df}

We now state the main result of this paper.
\begin{th}For weight transformation semigroup $(S,X,\omega)$ ; if
$S$ is topologically right simple and acts topologically
surjective on $X$, then $WAP(X,\omega) \subseteq LUC(X,\omega)$.\\
\textbf{Proof.}\\ \emph{ Let$((\epsilon ,\delta ) ,(T,Y,
\overline{\omega}))$ denotes the universal semitopological
compactification of $(S,X,\omega)$; where $\epsilon , \delta $ are
evaluation map on $S$ and $X$ respectively.\\
Put:~~~~~~~ $T_1=\{ t \in T : tT=T , tY=Y\}$. Then  $T_1$ is a
subsemigroup of $T$ contains $\epsilon(S)$ and by lemma [3.6]; for
each $f \in C(Y,\omega)$ the map $ _t(\frac{f}{\overline{\omega}})
~~_t{\overline{\Omega}} :T\longrightarrow C(Y)$ is norm continuous
at each point of
$\epsilon(S).$\\
For $s \in S $ and $f\in C(Y,\overline{\omega})$ we have:\\
$$\delta^\ast (
_{\epsilon(s)}(\frac{f}{\overline{\omega}})_{\epsilon(s)}\overline{\Omega})=
\hat L _s\delta^\ast(\frac{f}{\overline{\omega}})_s\Omega $$ where
$\delta^\ast :C(Y,\overline{\omega}) \longrightarrow
WAP(X,\omega)$ is the dual of evaluation map $\delta :X
\longrightarrow X^{WAP}.$ Therefore
$$WAP(X,\omega)=\delta^\ast (C(Y,\overline{\omega}))\subseteq
LUC(X,\omega).$$}
\end{th}
\begin{co} \emph{If $(S,X,\omega)$ is a weighted transformation
group, then :\\$WAP(X,\omega) \subseteq LUC(X,\omega)$}\\
 \textbf{Proof.}\\
 \emph{ Since $(S,X,\omega)$ is weighted transformation group
 by 1.1.17 of \cite{Berg}; $S$ is left simple and right simple.\\
 Let $e$ be the identity of $S$. For $s \in S$ and $x \in X$ we have:
 $ x=ex=ss^{-1}x \in SX$ and so $sX=X$. By theorem 3.8 the proof
 is obviously.}
\end{co}
The following example show that the condition of above theorem is
necessary.
\begin{ex}\emph{ If $S=(\mathbb{N},+), X=(Q,+)$, $\omega(x)=e^x$ and suppose that action of $S$ on $X$
 is $(s,x)\longrightarrow s+x$. Then $WAP(X,\omega)\varsubsetneqq
LUC(X,\omega)$. Since by 4.19 of \cite{Berg}, there exist $f_0 \in
WAP(X)$ such that $f_0 \notin LUC(X)$. Now assume that $f=f_0
\omega$. Therefore
 $f\in WAP(X,\omega)$ but $f \notin LUC(X,\omega).$}
\end{ex}

{\bf Acknowledgements}  
This research was supported by center of excellence for
mathematics at Isfahan university and first author would like to
thank from Razi University for his support.


\begin{thebibliography}{}
\bibitem{Berg} J. F. Berglund, D. H. Junghenn and P. Milnes; \emph{Analysis
on semigroups}, Canadian Math. Soc. Ser. Mono. Adv. Texts, Niley,
New York, 1989.
\bibitem{Burckel} R. B. Burckel; \emph{Weakly almost periodic functions on semigroups},
 Notes on mathematics and its applications, Gordon and Breach, New York
 (1970).
\bibitem{Dzinotyiweyi} H. A. M. Dzinotyiweyi; \emph{The
analogue of the group algebra for topological semigroups}, Pitman
Advanced Publishing Program, 1984.
\bibitem{vishki} H. R. Ebrahimi-vishki and M. A. Pourabdollah;
\emph{More on norm-continuity of weakly almost periodic
functions}, J. Sci. Univ. Tehran.\textbf{2} (1997) 29-33.
\bibitem{Granirer} E.E. Granirer; \emph{ Weakly almost periodic and
uniformly continuous funcionals on the fourier algebra of any
locally compact group}, Trans. Amer. Math. Soc. \textbf{189}
(1974) 371-382.
\bibitem{Hewitt} E. Hewitt and K. A. Ross; \emph{Abstract harmonic analysis}, vol. 1, Springer-verlag, Berlin (1983)
\bibitem{H} J. M. Howie; \emph{An introduction to semigroup theory},
Academic Press, London, (1976).
\bibitem{Junghenn} H. D. Junghenn; \emph{Almost periodic compactification of transformation semigroup}
, Pac. J. Math, \textbf{57 (1)} (1975) 207-216.
\bibitem{khadem} A. A. Khadem-Maboudi and M. A. Pourabdollah;
\emph{Almost periodic type function spaces on weighted
semitopological semigroup}, Farest J. Math. Sci. Special
Volume(1998)  part 1, 107-119.
\bibitem{Milnes} P. Milnes; \emph{Extention of continuous functions on topological semigroup}.
Pacific Journal of mathematics. \textbf{58 (2)} (1975), 553-562.
\bibitem{Milnes and Pym} P. Milnes, and J. Pym, \emph{Function spaces on semitopological
semigroups},
  Semigroup Forum \textbf{19} (1980), 347-354
\bibitem{Rejali} A. Rejali;\emph{Weighted function spaces on topological group},
Bull. Iranian math. soc. \textbf{22} (1998) 43-63.
\bibitem{Ruppert} W. Ruppert, \emph{Compact semitopological semigroups},An Intrinsic Theory, Lecture Notes in Mathematics.
\textbf{1079}, Springer-Verlag, New York (1984).\\\\





\footnotesize{The address of both authors:

Department of Mathematics, Isfahan University, Isfahan, Iran.

E-mail address: M.Abolghasemi@sci.ui.ac.ir

Department of Mathematics, Isfahan University, Isfahan, Iran.

E-mail address: Rejali@sci.ui.ac.ir

Faculty of Mathematics, Ferdowsi University of Mashhad, Mashhad
Iran.

E-mail address: Vishki@math.um.ac.ir}

\end{thebibliography}
\end{document}